\documentclass[12pt, reqno, a4paper]{amsart}

\usepackage{ amssymb, amsmath, enumerate, amsfonts, amsthm, mathrsfs, url, bm, mathtools, enumitem}

\usepackage{xcolor}  	
\usepackage{hyperref}
\hypersetup{
colorlinks,
   linkcolor={cyan!80!black},
   citecolor={cyan!80!black},
 urlcolor={cyan!80!black}
}

\usepackage{color}

\usepackage[margin=1in]{geometry}

\RequirePackage{doi}

\usepackage{amscd}
\usepackage{amsfonts}
\usepackage{float}
\usepackage{color}
\usepackage[
backend=biber,
style=alphabetic,
]{biblatex}
\usepackage{bookmark}

\renewbibmacro{in:}{}
\DeclareFieldFormat{title}{#1}

\DeclareFieldFormat[article]{title}{\mkbibemph{#1}}
\DeclareFieldFormat[incollection]{title}{\mkbibemph{#1}}
\DeclareFieldFormat[book]{title}{\mkbibemph{#1}}
\DeclareFieldFormat[incollection]{booktitle}{#1}
\DeclareFieldFormat[article]{journaltitle}{#1}

\AtEveryBibitem{\ifentrytype{misc}{\DeclareFieldFormat{title}{\mkbibemph{#1}}}{}}

\DeclareFieldFormat{eprint:eprint}{arXiv:\href{https://arxiv.org/abs/#1}{#1}}

\DeclareFieldFormat[inproceedings]{title}{\mkbibemph{#1}}

\DeclareFieldFormat[inproceedings]{booktitle}{#1}

\newtheorem{theorem}{Theorem}[section]
\newtheorem{lemma}{Lemma}[section]
\newtheorem{example}{Example}[section]
\newtheorem{corollary}{Corollary}[section]
\newtheorem{proposition}{Proposition}[section]

\theoremstyle{definition}
\newtheorem{definition}{Definition}[section]

\theoremstyle{remark}
\newtheorem{remark}{Remark}[section]

\numberwithin{equation}{section}

\newcommand{\Mod}[1]{\ (\mathrm{mod}\ #1)}

\renewcommand{\epsilon}{\varepsilon}

\renewcommand{\leq}{\leqslant}
\renewcommand{\geq}{\geqslant}

\renewcommand{\le}{\leqslant}
\renewcommand{\ge}{\geqslant}
\newcommand{\bi}[2]{\mathrm{BI}_{#1}\big[#2\big]}

\begin{document}

\title{Visiting early at prime times}

\author{Tony Haddad}
\address{Department of Mathematics and Statistics \\
University of Turku \\
20014 Turku \\
Finland}
\curraddr{}
\email{{\tt tohadd@utu.fi}}
\thanks{}

\author{Sun-Kai Leung}
\address{Mathematical Institute, University of Oxford, Andrew Wiles Building\\ Radcliffe Observatory Quarter, Woodstock Rd\\
Oxford OX2 6GG\\
United Kingdom}
\email{sunkaileung@gmail.com}

\author{Cihan Sabuncu}
\address{Max Planck Institut für Mathematik\\ 
Vivatsgasse 7 \\
53111 Bonn \\
Germany}
\curraddr{}
\email{sabuncu@mpim-bonn.mpg.de}
\thanks{}


\date{}


\begin{abstract}
Given an integer $m \geq 2$ and a sufficiently large $q$, we prove that there always exists an arithmetic progression with common difference $q$ for which the $m$-th least prime in such progression is $\ll q \exp(3.8075 m)$. This is obtained using a minorant for the indicator of primes in the Maynard--Tao sieve, in parallel with Baker--Irving and Stadlmann's work on bounded gaps between primes. The dependence on $q$ is best possible. Furthermore, we generalize our result to dynamical systems. The quality of the result depends crucially on the first return time, which we illustrate in the case of Diophantine approximation.

\end{abstract}

\maketitle

\section{Introduction}

Analytic Number Theory has witnessed many breakthroughs on fundamental problems in recent years, in particular the bounded gaps between primes, thanks to Zhang \cite{MR3171761}, Maynard \cite{MR3272929} and Tao (see also \cite{MR3373710}), building upon the celebrated work of Goldston--Pintz--Y\i ld\i r\i m \cite{MR2552109}. Maynard \cite{MR3530450} further proved that any subset of primes which is ``well distributed" in arithmetic progressions contains many primes that are close together. 

In particular, for primes in a fixed progression, the best known result is due to Baker and Zhao \cite{MR3632085}, who proved the following. Given $m \geq 1$, let $x$ be a sufficiently large real number, $q=x^{\theta}$ a suitably ``powerful" modulus with $0<\theta\leq 5/12-\eta$ for a fixed $\eta>0$ and $\prod_{p|q} p \leq (\log x)^C$ for a fixed $C>0$, and $(a,q)=1.$ Then there exist primes $p_1<\cdots<p_m$ in $(x/2,x]$ such that $p_i \equiv a \Mod{q}$ for $1 \leq i \leq m$ and
\begin{align*}
p_m-p_1 \leq q \exp(cm)
\end{align*}
for some constant $c=c(\theta)>0.$ 

Apart from progressions, there are also infinitely many 
bounded gaps between primes in a fixed ``Chebotarev set" \cite{MR3344173}, a suitable ``Sato--Tate interval" \cite{MR4048052} or a ``nil-Bohr set" \cite{MR4328791}.

To motivate our work, we begin by asking: how big is the first prime in a given arithmetic progression? Let $p(q, a)$ denote the least prime in the arithmetic progression $a \Mod{q}.$ Then recall Linnik's celebrated theorem that  $p(q,a) \ll q^L$ for some absolute constant $L>0.$  If $p_m(q, a)$ is now the $m$\textsuperscript{th} prime in the arithmetic progression $a \pmod q$, then a quantitative version of Linnik's theorem \cite[Corollary 18.8]{MR2061214} implies
\[
    \max_{\substack{a \Mod q\\(a,q)=1}} p_m(q,a)\ll_m q^L,
\]
for some $L>0$.

Assuming a uniform variant of the Hardy--Littlewood prime $k$-tuple conjecture, the second author \cite{moments} showed that as soon as $f(q) \to \infty$ with $q \to \infty,$ the least prime $p(q,a)$ is $\leq f(q)\varphi(q)\log q$ for almost all $a \Mod{q}$. In fact,
the least prime in arithmetic progressions to a common large modulus has an exponential limiting distribution. More precisely, if a reduced residue $a \Mod{q}$ is chosen uniformly at random, then we have the convergence in distribution
\begin{gather*}
\frac{p(q,a)}{\varphi(q)\log q} \xrightarrow[]{d} \operatorname{Exp}(1) \qquad \text{as $q \to \infty.$ }
\end{gather*}

Given $m \geq 1$ and $(a,q)=1,$ we are also interested in the $m$-th least prime $p \equiv a \pmod{q},$ which we denote as $p_m(q,a).$ For any residue class containing fewer than $m$ primes, we set $p_m(q,a)=\infty$. As $a \Mod{q}$ varies uniformly with $q \to \infty$, the second author \cite{MR4865836} further showed conditionally that the normalized primes in the progression, 
\begin{gather*}
\frac{p_1(q,a)}{\varphi(q)\log q} < \frac{p_2(q,a)}{\varphi(q)\log q} < \cdots,
\end{gather*}
are located on the positive real number line $(0,\infty)$ as if they were randomly positioned points, in the sense of a Poisson point process.

In this paper, we study $\min_{a \pmod q} p_m(q, a)$. By the pigeonhole principle, one can immediately show that $\min_{a \Mod{q}} p_m(q,a) \ll_m \varphi(q)\log q$. Applying a variant of the Maynard--Tao sieve weights, we improve this bound.

\begin{theorem} \label{thm:1.1}
For $m \geq 2$ and $q$ sufficiently large in terms of $m$, we have
\begin{align*}
\min_{a \Mod{q}} p_m(q,a) \ll q \cdot \exp(3.8075m).
\end{align*}
\end{theorem}

\begin{remark}
    A similar result can be obtained from the setting of Maynard \cite[Theorem 3.1]{MR3530450}, but this natural statement does not seem to have been observed before in the literature. Moreover, Maynard's theorem gives the upper bound $O(qe^{cm})$ for some large constant $c$. If we go into the proof and optimise the constant for our special case, we may take any $c>4$. As a by-product of our result, we see that we can take the exponent $c=3.8075$ like in Stadlmann \cite{Stadlmann2025}. 
\end{remark}

When keeping $m$ fixed and letting $q$ grow, the order of magnitude is optimal. Indeed, in any arithmetic progression, the $m$-th prime is at least the $m$-th positive integer in that progression. Thus $\min_{a \Mod q} p_m(q, a) > (m-1)q$.

The dependence on $m$ in Theorem \ref{thm:1.1} could be improved under the following variant of the Hardy--Littlewood conjecture: For all fixed $\epsilon > 0$ and fixed admissible $m$-tuple of minimal diameter $\mathcal H = \{0 = h_1 < \cdots < h_m\}$, we expect that for sufficiently large $q$, we must have 
\[
    \#\Big\{\epsilon q < a < 2\epsilon q : a+q h_i \text{ is prime for all $i$}\Big\} \sim \mathfrak S(\mathcal H, q)\frac{\epsilon q}{(\log q)^m}
\]
with
\[    
    \mathfrak S(\mathcal H, q) = \left(\frac{q}{\varphi(q)}\right)^{m-1}\prod_{p\nmid q} \left(1-\frac{\#\{h_i \Mod p\}}{p}\right)\left(1-\frac{1}{p}\right)^{-m}.
\]
This variant of the Hardy--Littlewood conjecture is consistent with the Cram\' er--Granville model (see \cite{MR1349149}). With this conjecture, we directly obtain 
\[
    \min_{a \pmod q} p_m(q, a) < (H(m)+o(1))q 
\]
as we let $q \to \infty$, with $H(k)$ being the diameter of the narrowest admissible $k$-tuple. Since $H(m) \asymp  m\log m$, we would have 
\begin{equation}
    \label{eq:conditional-bound}
    \min_{a \pmod q} p_m(q, a) \ll qm\log m,
\end{equation}
conditional on the variant of the Hardy-Littlewood conjecture.

The bound \eqref{eq:conditional-bound} would be sharp: Consider $a, q$ satisfying $(a, q) = 1$ and a prime modulus $q > m$. If $(a+\ell_iq)_{i=1}^m$ are the first $m$ primes in the arithmetic progression $a \pmod q$, then the set $\{\ell_2, \ldots, \ell_m\}$ is admissible. Therefore, taking the minimum over all integers $a$ coprime to $q$, we have 
\[
    \min_{a \pmod q} p_m(q, a) > (H(m-1)+1)q.
\]

Also, by varying the progressions instead of fixing a particular one, we are able to obtain a prime cluster early,\footnote{See \cite{MR0578815} and a recent improvement \cite{MR3698292} for a lower bound for the least prime in a progression instead.} which makes the nature of our result distinct from the existing literature, including those discussed in the beginning of the paper. 
    
Furthermore, we generalize our result to dynamical systems, for example circle rotations (Diophantine approximation, see Example \ref{eg:rot}). It turns out that the quality of the result depends crucially on the first return time, which we explore in Section \ref{sec:1streturn}. In particular, in the case of circle rotations $T_{\alpha}:x \mapsto x+\alpha$ for almost all $\alpha$ or those $\alpha$ with bounded partial quotients, we establish tighter lower and upper bounds than those in the literature for the first return time (see Proposition \ref{prop:1streturn}).

Recurrence is a central theme in ergodic theory, i.e., how points in dynamical systems return close to themselves under iteration. In 1890, Poincar\'{e} proved a simple yet far-reaching recurrence theorem in the context of the ``three-body" problem of planetary orbits (see \cite[p. 21]{MR2723325}). For our purposes, here we state a simple version (see \cite[Exercise 2.2.1]{MR2723325} for instance).

\begin{lemma}[Poincar\'{e} recurrence theorem] \label{lem:poincare}
Given a measure-preserving system $(X, \mathcal{B},\mu, T)$ with $\mu(X)=1$, let $A \in \mathcal{B}$
with $\mu(A) > 0$. Then there exists an integer $1 \leq n \leq \lfloor\mu(A)^{-1}\rfloor$ such that
\begin{align*}
\mu (A \cap T^{-n} A) >0.
\end{align*}
\end{lemma}

Indeed, if $\mu(A\cap T^{-n}A)=0$ for $1\leq n\leq N:=\lfloor\mu(A)^{-1}\rfloor$, then $A,T^{-1}A,\ldots,T^{-N}A$ are pairwise disjoint up to null sets, contradicting $(N+1)\mu(A)>1$. A famous example is Dirichlet's approximation theorem, corresponding to circle rotation. Inspired by the recurrence theorem, we shall visit early at prime times in \textit{metric-measure-preserving systems} defined as follows.

\begin{definition}[Metric-measure-preserving system] \label{def:mmps}
A \textit{metric-measure-preserving system} is a quintuple $(X,d,\mathcal{B},\mu,T)$ where $(X,d)$ is a metric space and $(X,\mathcal{B},\mu,T)$ is a measure-preserving system satisfying the following compatibility conditions:
\begin{itemize}
    \item the Borel $\sigma$-algebra $\mathcal{B}$ is induced by the metric $d;$ 
    \item the probability measure $\mu$ is \textit{pseudo-doubling}, i.e., there exist constants $\lambda \geq 1$ and $r_0>0$ such that
    \begin{align*}
     \mu(B(x; 2r)) \leq \lambda \cdot \mu(B(x; r))
    \end{align*}
    for any $x \in X$ and $0<r<r_0$, where $B(x;r):=\{ y \in X\,:\, d(y,x) <r \}$;
    \item the measure-preserving map $T$ is an isometry.
\end{itemize}
\end{definition}

Let $p_m^T(x_0,x;\varepsilon)$ denote the $m$-th least prime $p$ for which $T^p x_0 \in B(x;\varepsilon).$ If fewer than $m$ primes satisfy this condition, we set $p_m^T(x_0,x;\varepsilon)=\infty$. Then, with the definition in hand, we state the dynamical generalization of Theorem \ref{thm:1.1}.

\begin{theorem} \label{thm:1.2}
Let $m \geq 2$ and $(X,d,\mathcal{B},\mu,T)$ be a metric-measure-preserving system. Then for $0<\varepsilon<r_0$ sufficiently small (in terms of $m$) and $x_0 \in \operatorname{supp}\mu$, we have
\begin{align*}
 \min_{x \in X} \, p_m^T(x_0,x;\varepsilon) \ll_{m,\lambda} {\mu(B(x_0;\varepsilon))}^{-1}.
\end{align*}
\end{theorem}

\noindent\textit{Notation.} Throughout the paper, we use the standard big $O$ and little $o$ notations as well as the Vinogradov notations $\ll, \gg$ and the Hardy notation $\asymp,$ where the implied constants depend only on the subscripted parameters.

\section{Main results}

We start by defining families of minorants for the indicator function of prime numbers. This definition is inspired by a similar one given by Baker and Irving in \cite{BakerIrving2017}.
For $n\geq2$, let $P^-(n)$ denote the least prime factor of $n$, and set $P^-(1)=\infty$.
\begin{definition}[Minorants for the indicator of primes]
    \label{def:set of minorants}
    Fix a positive integer $\ell$, parameters $(\theta, \xi, \beta) \in (0, 1]^3$, and a function $\delta \colon (0, \theta) \to (0, \xi]$. We define $\bi{\ell}{\theta, \xi, \beta, \delta}$ as the set of functions $\rho: \mathbb N \times (1, \infty) \to \mathbb R$ satisfying
    \begin{enumerate}[label={\textnormal{(\roman*)}}]
    \item \label{item:minorant for primes} There exists $C > 0$ such that $\rho(a, q) \le 1_\mathbb P(a)$ for every $a > q > 1$, and $-C \le \rho(a, q)$ whenever $q > 1$ and $\ell q < a \le (\ell+1)q$.

    \item \label{item:supported on rough} For every $a, q$ satisfying $a > q > 1$ and $\rho(a, q) \ne 0$, we have $P^-(a) \ge q^\xi$.

    \item \label{item:level of distribution} For every $q \ge 2$, $c \in \mathbb Z$, $A > 0$ and $\epsilon \in (0, \theta)$, we have 
    \[
        \sum_{\substack{r \le q^{\theta-\varepsilon} \\ r|P}} \Bigg|\sum_{\ell q < a \le (\ell+1) q} \rho(a, q)\Big(1_{a \equiv c \pmod r} - \frac{1_{(a, r) = 1}}{\varphi(r)}\Big)\Bigg| \ll_{A, \varepsilon} \frac{q}{(\log q)^A}
    \]
    with $P = \prod_{\substack{p\le q^{\delta(\epsilon)} \\ p\nmid c}} p$.

    \item \label{item:pnt for minorant}We have 
    \[
        \sum_{\ell q < a \le (\ell+1) q} \rho(a, q) \ge (\beta - o(1))\frac{q}{\log q}
    \]
    as $q \to \infty$.
\end{enumerate}
\end{definition}

Baker and Irving \cite[Lemma 1]{BakerIrving2017} proved the following theorem on bounded gaps between primes: If $\bi{1}{\theta, \xi, \beta, \delta}$ is not empty for some parameters $(\theta, \xi, \beta) \in (0, 1]^3$ and some function $\delta \colon (0, \theta) \to (0, \xi]$, then 
\begin{equation}
\label{eq:bakerirving}
    H_m \coloneqq \liminf_{n \to \infty} (p_{n+m} - p_n) \ll_\epsilon \exp\bigg(\frac{(2+\epsilon)m}{\theta\beta}\bigg)
\end{equation}
where $p_1 < p_2 < \cdots$ is the increasing sequence of prime numbers.

Different theorems in the literature are proofs that the indicator $1_\mathbb P$ is in $\bi{1}{\theta, 1, 1, \delta}$ for some $\theta$ and $\delta$, leading to proofs of the form $H_m \ll_{\epsilon} e^{(2/\theta+\epsilon)m}$: 
\begin{itemize}
    \item The Bombieri-Vinogradov Theorem implies that we can take $\theta = \frac{1}{2}$ and the function $\delta \equiv 1$. We note that Maynard \cite{MR3272929} proved that $H_m \ll m^3e^{4m}$. The exponent $4$ comes from applying the Bombieri-Vinogradov Theorem in a key part of the proof.

    \item Zhang \cite{MR3171761} essentially proved that we can take $\theta = \frac{1}{2} + \frac{1}{584}$ and any function $\delta$ satisfying $\delta(\epsilon) \in (0, \frac{1}{1168})$.

    \item Polymath \cite{polymath-zhangtype} proved that we can take $\theta = \frac{1}{2} + \frac{7}{300}$ and any $\delta$ satisfying $\delta(\epsilon) < \frac{5\epsilon}{3}$.

    \item Stadlmann \cite{Stadlmann2025} proved that we can take $\theta = \frac{1}{2} + \frac{1}{40}$ and any $\delta$ satisfying $\delta(\epsilon) < \frac{3\epsilon}{2}$.
\end{itemize}

With the Elliott-Halberstam conjecture, we expect that $1_\mathbb P \in \bigcap_{\ell \ge 1} \bi{\ell}{1, 1, 1, 1}$. (It is important to note here that when taking the intersection, the implicit constant in Definition \ref{def:set of minorants} \ref{item:level of distribution}, as well as the speed of convergence for Definition \ref{def:set of minorants} \ref{item:pnt for minorant} are not necessarily the same for all $\ell \ge 1$.)

Concerning non-trivial minorants of $1_\mathbb P$, we have that 
\begin{itemize}
    \item Baker and Irving \cite{BakerIrving2017} showed that for all $\epsilon > 0$, there exist $\xi$ and $\delta$ such that $\bi{1}{\frac{1}{2} + \frac{16}{659} - \epsilon, \xi, 1-8\cdot 10^{-6}, \delta}$ is not empty.

    \item Stadlmann \cite{Stadlmann2025} constructed the following minorant:
    \begin{equation}\label{eq:stadlmann minorant}
    \rho_{0}(a, q) \coloneqq 1_\mathbb P(a) 
    - \sum_{\substack{a = p_1p_2p_3p_4n_5 \\ 0.19038 \le \alpha_4 < \alpha_3 < \alpha_2 < \alpha_1 \le 0.40481 \\ \alpha_1+\alpha_2 < 0.40481 \\ \alpha_2+\alpha_3+\alpha_4 > 0.59519 \\ P^-(n_5) \ge p_4}} 1 -  
    \sum_{\substack{a = p_2p_3p_4p_5p_6 \\ 0.19038 \le \alpha_2, \alpha_3, \alpha_4, \alpha_5, \alpha_6 \le 0.23848 \\ \alpha_2, \alpha_4 > \alpha_3 \\ \alpha_2+\alpha_4 < 0.40481 \\ \alpha_2+\alpha_3+\alpha_5 > 0.59519 \\ \alpha_6 \ge \alpha_5}} 1,
\end{equation}
where $\alpha_j = \frac{\log p_j}{\log q}$ in the sums above. She proved that for $\delta \equiv 10^{-10}$, we have $\rho_0 \in \bi{1}{0.5253, 0.19038, 0.99998, \delta}$.
\end{itemize}

We can now state an analogous statement to Baker and Irving's result \eqref{eq:bakerirving} for primes in arithmetic progressions, which provides a lower bound for the number of residue classes with early prime clusters.

\begin{theorem} \label{thm:main}
    Fix $\epsilon > 0$ and let $m \geq 2$. Suppose that $\bigcap_{\ell \ge 1} \bi{\ell}{\theta, \xi, \beta, \delta}$ is not empty for some parameters $(\theta, \xi, \beta) \in (0, 1]^3$ and some function $\delta \colon (0, \theta) \to (0, \xi]$ satisfying $\delta(y) \gg y$. Then for any sufficiently large modulus $q$ (in terms of $m$ and $\epsilon$), there exist constants $C, C'>0$ which depend on $(\theta, \xi, \beta, \delta, \epsilon)$ such that
    \begin{align*}
        \# \bigg\{ a \Mod{q} \,: \, p_m(q,a) \leq C q \exp\Big(\tfrac{2}{\theta\beta}(1+\epsilon) m\Big) \bigg\}
        \gg_{m,\epsilon} \varphi(q)(\log q)^{-C' \exp( \tfrac{2}{\theta\beta}(1+\epsilon) m)}.
    \end{align*}
\end{theorem}

Now, we can state our theorem in the case of Stadlmann's minorant. For this, we have the following lemma.
\begin{lemma}
    Let $\rho_0$ be defined as in \eqref{eq:stadlmann minorant}. For the constant function $\delta(\epsilon) = 10^{-10}$, we have 
    \[
        \rho_0 \in \bigcap_{\ell \ge 1} \bi{\ell}{0.5253, 0.19038, 0.99998, \delta}.
    \]
\end{lemma}

\begin{proof}
    Properties \ref{item:minorant for primes} and \ref{item:supported on rough} of Definition \ref{def:set of minorants} follow directly from the definition of $\rho_0$. Concerning properties \ref{item:level of distribution} and \ref{item:pnt for minorant}, Stadlmann proved in \cite[Section 4]{Stadlmann2025} that $\rho_0 \in \bi{\ell}{0.5253, 0.19038, 0.99998, \delta}$ for $\ell = 1$. For each fixed $\ell \ge 2$, we apply the same decomposition and distribution estimates on $(\ell q,(\ell+1)q]$. Since $\log(\ell q)=\log q+O_\ell(1)$, replacing $\log q$ by $\log(\ell q)$ changes each exponent $\alpha_j=\log p_j/\log q$ by $O_\ell(1/\log q)$. For fixed $\ell$, these $O_\ell(1/\log q)$ perturbations are eventually smaller than the fixed slack in the Type I/II/III ranges used in Stadlmann's argument, so the same distribution estimates remain applicable after slightly enlarging the corresponding ranges. Their constants, thresholds, and $o(1)$ rates may depend on $\ell$, while the numerical parameters $0.5253$, $0.19038$, $0.99998$ and $\delta\equiv10^{-10}$ remain unchanged.
\end{proof}
Using this, we get
\begin{corollary}\label{stadlman_corollary}
We have for $q$ sufficiently large in terms of $m$,
$$ \# \big\{ a \Mod{q} \,: \, p_m(q,a) \leq C q \exp(3.8075 m) \big\}
        \gg_{m} \varphi(q)(\log q)^{-C' \exp( 3.8075 m)}. $$
\end{corollary}

We also state the quantitative version of Theorem \ref{thm:1.2}, which is
the dynamical generalization of Theorem \ref{thm:main}.

\begin{theorem} \label{thm:dynamic}
Let $m \geq 2$ and $(X,d,\mathcal{B},\mu,T)$ be a metric-measure-preserving system. Then 
there exist constants $C=C(\lambda), C'=C'(\lambda)>0$ such that for any nonperiodic $x_0 \in \operatorname{supp}\mu$ and $0<\varepsilon<r_0$ sufficiently small (in terms of $m$ and $x_0$), we have
\begin{gather*}
\mu ( \{ x \in X \,: \, p_m^T(x_0,x;\varepsilon) \leq C q  \exp( 3.8075 m) \} ) \\
\gg_{m,\lambda} \varphi(q) (1+\log q)^{-C'\exp(3.8075m)} \mu(B(x_0;\varepsilon)) 
\end{gather*}
for some integer $1 \leq q \ll_\lambda m^{\log_2 \lambda} \exp( 3.8075 \cdot  m \log_2 \lambda ) {\mu(B(x_0;\varepsilon))}^{-1}$.

\end{theorem}

\begin{remark}
    In the proof of Theorem \ref{thm:dynamic}, $q$ is the first return time of $x_0$ to $B(x_0;\varepsilon/(2h))$, with $h$ defined there.
\end{remark} 
\begin{remark}
We could also state our theorem in the style of Theorem \ref{thm:main} using the minorants, but we omit this for the sake of clarity.
\end{remark}

We provide an example to illustrate this theorem.

\begin{example}[Circle rotation] \label{eg:rot}
Consider the metric-measure-preserving system with $X = \mathbb R/\mathbb Z, d=\|\cdot \|, \mu$ being the normalized Haar measure and $T_{\alpha}:x \mapsto x+\alpha \Mod{1},$ where $\| \cdot \|$ is the distance to the nearest integer function and $\alpha$ irrational in $[0,1].$
Let $m \geq 2.$ Then for any sufficiently small $\varepsilon>0$ (in terms of $m$), we have
\begin{align*}
\min_{\beta \in [0,1]} p_m^{\alpha}(\beta;\varepsilon) \ll \varepsilon^{-1} m \exp(7.615m)
\end{align*}
where $p_m^{\alpha}(\beta;\varepsilon)$ denotes the $m$-th least prime $p$ for which $\| p \alpha  -\beta\|<\varepsilon.$ If fewer than $m$ primes satisfy this condition, we set $p_m^{\alpha}(\beta;\varepsilon)=\infty$.
 Moreover, for a typical irrational $\alpha \in [0,1]$ and $\varepsilon>0$ sufficiently small in terms of $m$ and $\alpha$, we can establish a more precise version of Theorem \ref{thm:dynamic} in this case by showing that the first return time is $\asymp \varepsilon^{-1+o(1)}$; see Section \ref{sec:1streturn} and Proposition \ref{prop:1streturn} for details.

\end{example}

\section{A variant of Maynard--Tao sieve weight} 
In Sections 3--5, $\mu$ denotes the Möbius function; in the dynamical sections, it denotes the probability measure.

Let $m, k \geq 2$. Let $(h_1,\ldots,h_k)$ be an admissible $k$-tuple, i.e., a $k$-tuple of non-negative integers satisfying $0 \le h_1<\cdots<h_k$ and $\# \{ h_i \Mod{p} \,:\, 1 \leq i \leq k \}<p$ for any prime $p$. Fix also $\theta, \xi, \beta$ and $\delta$ such that $\bigcap_{\ell \ge 1} \bi{\ell}{\theta, \xi, \beta, \delta}$ is not empty.

Then, we would like to construct a sequence of non-negative weights $w_a$ over the reduced residue system $\{a \in (q, 2q] \,:\, (a, q)=1 \}$ for any $q=q(m)$ sufficiently large, so that
\begin{align} \label{eq:pigeon}
S:= \sum_{\substack{q < a \le 2q  \\ (a,q)=1 }} \Bigg( \sum_{i=1}^k \rho (a+qh_i, q) - (m-1) - k \sum_{i=1}^k \sum_{\substack{p| a+qh_i \\ p \leq q^{\gamma}  }} 1 \Bigg) w_a
\end{align}
is positive for some small $\gamma=\gamma(k)>0$ and a minorant $\rho \in \bigcap_{\ell \ge 1} \bi{\ell}{\theta, \xi, \beta, \delta}$. Let $\Omega(n)$ count the prime factors of $n$ with multiplicity. We define
\begin{align}
\begin{split}
    \label{eq:def of Mz}
    y_a &:= \sum_{i=1}^k 1_\mathbb P(a+qh_i) - (m-1) - k \sum_{i=1}^k \sum_{\substack{p| a+qh_i \\ p \leq q^{\gamma}  }} 1 \\
    M(z) &:= \max\Big\{w_a : \text{$a \in (q, 2q]$ satisfying $\Omega(a+qh_i) < z$ for all $i$}\Big\}
\end{split}
\end{align}
Suppose that $q \ge 2+h_k$. If $y_a \ge 1$, then the finite sequence $(a+qh_i)_{i=1}^k$ must contain at least $m$ primes and all terms of the sequence must satisfy $P^-(a+qh_i) > q^\gamma$. It therefore implies that $\Omega(a+qh_i) < 2/\gamma$ for all $i$. Thus, we obtain the following lower bound
\begin{multline}\label{eq:count}
    \#\Big\{a \in (q, 2q] : \text{at least $m$ of the $a+qh_i$ are primes}\Big\} \\ \ge \sum_{\substack{q < a \le 2q \\ (a, q) = 1 \\ y_a \ge 1}} 1 
    \ge \sum_{\substack{q < a \le 2q \\ (a, q) = 1 \\ y_a \ge 1}} \Big(\frac{y_a}{k}\Big) \cdot \Big(\frac{w_a}{M(2/\gamma)}\Big) \ge \frac{S}{k\cdot M(2/\gamma)}.
\end{multline}

To construct our weight $w_a,$ let $w:=\log\!\log\!\log q$ and $W_q:=\prod_{\substack{p \leq w \\p \nmid q}} p.$ By the admissibility of the $k$-tuple $(h_1,\ldots,h_k),$ we can choose $b_0 \Mod{W_q}$ so that for those $a \equiv b_0 \Mod{W_q},$ we have $a+qh_i \equiv b_0+qh_i\not\equiv 0 \Mod{p}$ for any prime $p \leq  w, p \nmid q.$ Then, analogously to the Maynard--Tao sieve weight (see \cite[p.11]{MR3373710}), we let
\begin{align}
\label{eq:w_a def}
w_a:=1_{a \equiv b_0 \Mod{W_q}} \Bigg( \underset{\substack{d_i | a+qh_i\\i=1,\ldots,k}} {\sum \cdots \sum}
\prod_{i=1}^k \mu(d_i)  \times F\left( \frac{\log d_1}{\log q} ,\cdots, \frac{\log d_k}{\log q}  \right) \Bigg)^2
\end{align}
for some function $F:[0,\infty)^{k} \to \mathbb{R}$ supported on the truncated simplex
\begin{align*}
\Delta_k(\theta, \delta):=\Big\{(t_1,\ldots,t_k) \in \big[0, \delta(\tfrac{\theta}{4\log k})\big]^k \,:\, t_1+\cdots+t_k \leq \tfrac{\theta}{2}\big(1-\tfrac{1}{2\log k}\big) \Big\}.
\end{align*}
We require $F$ to be a finite linear combination of smooth tensor products $\prod_{i=1}^k f_i(t_i)$, with $f_i\in\mathcal C_c^\infty([0,\infty))$, such that each tensor product has support contained in the interior of $\Delta_k(\theta,\delta)$ taken within $[0,\infty)^k$. For future reference, the family of such functions is denoted as $\mathcal{F}_{k}(\theta, \delta)$.

With our choice of weight $w_a,$ it turns out that not only is the sum 
$S$ in (\ref{eq:pigeon}) computable, but it is also positively large.

To simplify notation, given any square-integrable function $G:[0,\infty)^k \to \mathbb{R},$ let us denote the singular integrals
\begin{align*}
\mathfrak{I}(G):=\int_{[0,\infty)^k} G(t_1,\ldots,t_k)^2 dt_1 \cdots dt_k
\end{align*}
and
\begin{align*}
\mathfrak{J}_i(G):=\int_{[0,\infty)^{k-1}} \left(\int_{[0,\infty)}G(t_1,\ldots,t_k) dt_i \right)^2 dt_1 \cdots dt_{i-1} dt_{i+1} \cdots dt_k
\end{align*}
for $1 \leq i \leq k.$ Then
in view of (\ref{eq:pigeon}), we are required to estimate the following sums.

\begin{proposition}[Non-prime sum] \label{prop:nonprime}
Fix $k \ge 2$, $\theta \in (0, 1]$, $\delta \colon (0, \theta) \to (0, 1]$, and $F \in \mathcal F_{k}(\theta, \delta)$. We define $w_a$ as in \eqref{eq:w_a def}. Then as we let $q \to \infty$, we have
\begin{gather*}
\sum_{\substack{q < a \le 2q\\(a,q)=1}} w_a
=\big(\mathfrak{I}(\partial_{t_1\cdots t_k}F)+o(1)\big)\frac{\varphi(q)}{W_q}\left( \frac{\varphi(W_q)}{W_q} \times  \frac{\varphi(q)}{q}  \times \log q\right)^{-k}.
\end{gather*}
\end{proposition}

\begin{proposition}[Prime sum] \label{prop:prime}
Fix $k \ge 2$, $1 \leq i_0 \leq k$, $(\theta, \xi, \beta) \in (0, 1]^3$, $\delta \colon (0, \theta) \to (0, \xi]$ and $F \in \mathcal F_{k}(\theta, \delta)$. We define $w_a$ as in \eqref{eq:w_a def}. Suppose that the set $\bigcap_{\ell \ge 1} \bi{\ell}{\theta, \xi, \beta, \delta}$ is not empty, and fix $\rho$, an element of this set. Then as $q \to \infty$, we have
\begin{gather*}
    \sum_{\substack{q < a \le 2q\\(a,q)=1}} \rho(a+qh_{i_0}, q) w_a 
    \ge \big( \beta \cdot \mathfrak{J}_{i_0}(\partial_{t_1\cdots t_k}F)-o(1)\big) \frac{\varphi(q)}{W_q}\left( \frac{\varphi(W_q)}{W_q} \times  \frac{\varphi(q)}{q}  \times \log q\right)^{-k}.
\end{gather*}

\end{proposition}

\begin{proposition}[Small prime factor sum] \label{prop:small_prime}
Fix $k \ge 2$, $1 \le i_0 \le k$, and $F \in \mathcal F_{k}(\theta, \delta)$ for some $\theta \in (0, 1]$ and $\delta \colon (0, \theta) \to (0, 1]$. Let $0<\gamma \leq \frac{1}{100k}$ and $q \ge 2$. Then 
\begin{align*}
    \sum_{\substack{q < a \le 2q \\ (a,q)=1 }} \bigg(\sum_{\substack{p| a+qh_{i_0}  \\ p \leq q^{\gamma} }} 1 \bigg) w_a \ll_{F} \gamma \times \frac{\varphi(q)}{W_q} \Bigg( \frac{\varphi(W_q)}{W_q} \times \frac{\varphi(q)}{q} \times \log q \Bigg)^{-k}.
\end{align*}
\end{proposition}

\section{Proof of Theorem \ref{thm:main}}

In this section, we shall prove Theorem \ref{thm:main} assuming the aforementioned propositions. All the implicit constants in this section are allowed to depend on $F, \theta, \xi, \beta$ and $\delta$. 

We start by fixing $(\theta, \xi, \beta) \in (0, 1]^3$ and $\delta \colon (0, \theta) \to (0, \xi]$ such that the intersection $\bigcap_{\ell \ge 1} \bi{\ell}{\theta, \xi, \beta, \delta}$ is non-empty, and we choose $\rho$ as an element of this intersection. Put
\[
R_k:=\frac{\theta}{2}\left(1-\frac{1}{2\log k}\right).
\]
Since $\delta(y)\gg y$, choose $d>0$ such that $\delta(y)\ge dy$ for all sufficiently small $y$, and fix $0<b<d/2$. The argument proving \cite[Theorem 40]{MR3373710} gives, for all sufficiently large $k$, a nonzero function $H_k$ supported in
\[
\big\{(u_1,\ldots,u_k)\in[0,b/\log k]^k:u_1+\cdots+u_k\le 1\big\}
\]
such that $\sum_i\mathfrak J_i(H_k)/\mathfrak I(H_k)\ge\log k-O_b(1)$. Indeed, in the notation of that proof, all we need is to take $T=b/\log k$, $\tau=1/\log k$, and $c=1/\log k+\alpha/\log^2 k$ with a fixed sufficiently negative $\alpha$. Rescale by $R_k(1-1/\log^2 k)$. The resulting support satisfies $\sum_i t_i\le R_k(1-1/\log^2 k)<R_k$ and $t_i\le bR_k(1-1/\log^2 k)/\log k<\delta(\theta/(4\log k))$ for every $i$. We then adapt \cite[proof of Theorem 22, (72)--(84)]{MR3373710} to approximate the rescaled density in $L^2$ by a finite sum of smooth tensor products. Each product can be supported in a rectangle $\prod_i[a_i,b_i]$ with $\sum_i b_i<R_k$ and $b_i<\delta(\theta/(4\log k))$ for every $i$. Since the ratio is continuous under $L^2$ approximation on this bounded region, we obtain a density $G_k$ with $\mathfrak I(G_k)=1$ and
\begin{equation}
\label{eq:LB for ratio}
\sum_{i=1}^k\mathfrak J_i(G_k)\ge\frac{\theta}{2}\log k-C_2
\end{equation}
for some $C_2>0$ depending only on $\theta$ and $\delta$. Define
\begin{equation}
\label{eq:def of F}
F(t_1,\ldots,t_k):=\int_{t_1}^{\infty}\cdots\int_{t_k}^{\infty}G_k(u_1,\ldots,u_k)\,du_k\cdots du_1.
\end{equation}
Then $F\in\mathcal F_k(\theta,\delta)$ and $\partial_{t_1\cdots t_k}F=(-1)^kG_k$.

Let $k_\ast$ be large enough for \eqref{eq:LB for ratio} to hold, and choose
\begin{equation}
\label{eq:choice for k}
k=\max\left\{k_\ast,\left\lceil\exp\left(\frac{2}{\theta\beta}(1+\epsilon/2)m+\frac{2C_2}{\theta}\right)\right\rceil\right\}.
\end{equation}
For this fixed $k$ and $F$, put
\[
D:=\beta\sum_{i=1}^k\mathfrak J_i(G_k)-(m-1)\mathfrak I(G_k)\ge 1+\epsilon m/2>0.
\]
Fix $C_F>0$ for the bound in Proposition \ref{prop:small_prime}, and choose
\[
0<\gamma\le\min\left\{\frac{1}{100k},\frac{D}{4C_Fk^2}\right\}.
\]
Define $S$ as in \eqref{eq:pigeon} with these choices. Combining Propositions \ref{prop:nonprime}, \ref{prop:prime} and \ref{prop:small_prime} gives
\begin{gather}
\label{eq:finalS}
S\ge\left(D-C_Fk^2\gamma-o(1)\right)
\frac{\varphi(q)}{W_q}\left(\frac{\varphi(W_q)}{W_q}\times\frac{\varphi(q)}{q}\times\log q\right)^{-k}.
\end{gather}
Taking $q$ large enough that $o(1)\le D/4$, we obtain
\begin{equation}
\label{eq:lowerS}
S\ge\frac D2\frac{\varphi(q)}{W_q}\left(\frac{\varphi(W_q)}{W_q}\times\frac{\varphi(q)}{q}\times\log q\right)^{-k}.
\end{equation}

Finally, we need to bound the weights $w_a$ over $\Omega(a+h_i q ) < 2/\gamma$ for all $i $. If $a \in (q, 2q]$ satisfies the bound on the number of prime factors $\Omega(a+qh_i) < z$ for all $i$, then
\[
    \frac{w_a}{\max_{(t_1, \ldots, t_k) \in \Delta_k(\theta, \delta)} F(t_1, \ldots, t_k)^2 } \le \bigg(\prod_{i=1}^k\sum_{d|a+qh_i}\mu^2(d)\bigg)^{2}  \le 4^{zk}.
\]
Since $F$ is smooth and compactly supported, $\|F\|_\infty<\infty$. It follows that $M(z)\leq\|F\|_\infty^2 4^{kz}$ with $M(z)$ defined as in \eqref{eq:def of Mz}. Combining this with \eqref{eq:count} and \eqref{eq:lowerS}, we get
\begin{equation*}
    \#\big\{ a \Mod{q} : \text{at least $m$ of the $a+qh_i$ are primes}  \big\} \gg_{m,\epsilon} \frac{\varphi(q)}{W_q}\left( \frac{\varphi(W_q)}{W_q} \times  \frac{\varphi(q)}{q}  \times \log q\right)^{-k}.
\end{equation*}

Finally, if $(a+qh_i)_{i=1}^k$ has at least $m$ primes, it follows that $p_m(q, a) \le q(2+h_k)$. As explained in \cite{MR3272929}, if we select $h_1, \ldots, h_k$ to be the first $k$ primes that are larger than $k$, then $(h_1, \ldots, h_k)$ is an admissible $k$-tuple that satisfies $h_k \ll k\log k$. The theorem follows.

\section{Proof of Propositions \ref{prop:nonprime}, \ref{prop:prime} \& \ref{prop:small_prime}}

By the definition of $\mathcal F_k(\theta,\delta)$, $F$ is a linear combination of tensor products, i.e. there are constants $c_j \in \mathbb{R}$ and $f_{i,j}\in \mathcal{C}^{\infty}_{c}([0,\infty))$ for $1 \leq i \leq k$ and $1 \leq j \leq J$ such that
\begin{equation}
\label{eq:F as a tensor prod}
F(t_1,\ldots, t_k)=\sum_{j=1}^J c_j\prod_{i=1}^k f_{i,j}(t_i).
\end{equation}

Let us denote
\begin{align*}
\lambda_{f}(n):=\sum_{d \mid n} \mu(d)f\left( \frac{\log d}{\log q} \right).
\end{align*}
Then our weight becomes
\begin{align*}
w_a &= 1_{a \equiv b_0 \Mod{W_q}} 
\left( \sum_{j=1}^J c_j \prod_{i=1}^k \lambda_{f_{i,j}} (a+qh_i) \right)^2 \\
&= 1_{a \equiv b_0 \Mod{W_q}} \sum_{j=1}^J \sum_{j'=1}^J c_j c_{j'} \prod_{i=1}^k
\lambda_{f_{i,j}}(a+qh_i)\lambda_{f_{i,j'}}(a+qh_{i}).
\end{align*}

\begin{proof}[Proof of Proposition \ref{prop:nonprime}]
By linearity, it suffices to estimate
\begin{align} 
T \coloneqq \sum_{\substack{q < a \le 2q\\(a,q)=1}} 1_{a \equiv b_0 \Mod{W_q}} \prod_{i=1}^k 
\lambda_{f_{i}}(a+qh_i) \lambda_{g_i}(a+qh_i)
\end{align}
for all $f_{i}, g_i \in \mathcal{C}^{\infty}_{c}([0,\infty))$ with $1 \leq i \leq k$, with both $\prod_{i=1}^k f_i(t_i)$ and $\prod_{i=1}^k g_i(t_i)$ being supported in $\Delta_k(\theta, \delta)$. After interchanging the order of summation, this becomes
\begin{align}\label{eq:prop3.1}
T = \underset{d_1, d_1' \ldots, d_k, d_k' \geq 1 }{\sum \cdots \sum} \left( \prod_{i=1}^k 
\mu(d_i) \mu(d_i') f_i \left( \frac{\log d_i}{\log q} \right)
g_i \left( \frac{\log d_i'}{\log q} \right) \right) S(\boldsymbol{d}, \boldsymbol{d'}),
\end{align}
where $\boldsymbol{d}:=(d_1,\ldots,d_k), \boldsymbol{d'}:=(d_1',\ldots,d_k')$ and
\begin{align*}
S(\boldsymbol{d}, \boldsymbol{d'}) \coloneqq
\sum_{\substack{q < a \le 2q\\ [d_i,d_i'] \, | \, a+qh_i, 1 \leq i \leq k\\a \equiv b_0 \Mod{W_q}}}
1_{(a, q) = 1}.
\end{align*}

Since $(h_1,\ldots,h_k)$ is fixed and $w=\log\!\log\!\log q$, we take $q$ sufficiently large that $w>\max_{i\ne j}|h_i-h_j|$.
Note that the sum $S(\boldsymbol{d}, \boldsymbol{d'})$ vanishes unless $q, W_q, [d_1,d_1'], \ldots, [d_k, d_k']$ are pairwise coprime. In this case, let $r_{W_q,\boldsymbol{d}, \boldsymbol{d'}}:=W_q\prod_{i=1}^k [d_i, d_i']$. It follows from the Chinese remainder theorem that there exists $c_{W_q, \boldsymbol{d}, \boldsymbol{d'}}$ such that
\[
S(\boldsymbol{d}, \boldsymbol{d'}) = \sum_{\substack{q < a \le 2q \\ a \equiv c_{W_q, \boldsymbol{d}, \boldsymbol{d'}} \Mod{r_{W_q, \boldsymbol{d}, \boldsymbol{d'}}}}} 1_{(a, q) = 1} 
\]
By Möbius inversion, we arrive at
\[
    S(\boldsymbol{d}, \boldsymbol{d'}) = \sum_{e|q}\mu(e) \sum_{\substack{q < a \le 2q \\ a \equiv c_{W_q, \boldsymbol{d}, \boldsymbol{d'}} \Mod{r_{W_q, \boldsymbol{d}, \boldsymbol{d'}}}}} 1_{e|a}
\]
From our previous discussion, we know that $(q, r_{W_q, \boldsymbol{d}, \boldsymbol{d'}}) = 1$. Hence,
\begin{equation}
\label{eq:estimate of S(d,d')}
    S(\boldsymbol{d}, \boldsymbol{d'}) = \sum_{e|q} \mu(e)\left(\frac{q}{er_{W_q, \boldsymbol{d}, \boldsymbol{d'}}} + O(1)\right) = \frac{\varphi(q)}{r_{W_q, \boldsymbol{d}, \boldsymbol{d'}}} + O(\tau(q)).
\end{equation}
Therefore, we define $T_1$ (the main term) and $T_2$ (the error term) as 
\begin{align*}
T_1 \coloneqq \frac{\varphi(q)}{W_q} \underset{\substack{d_1, d_1' \ldots, d_k, d_k' \geq 1\\q, W_q, [d_i,d_i'] \text{ pairwise coprime}}}{\sum \cdots \sum} \prod_{i=1}^k 
\frac{\mu(d_i) \mu(d_i')}{[d_i,d_i']} f_i \left( \frac{\log d_i}{\log q} \right)
g_i \left( \frac{\log d_i'}{\log q} \right),
\end{align*}
and $T_2 \coloneqq T - T_1$.

We rewrite $T_1$ as
\begin{align*}
T_1 = \frac{\varphi(q)}{W_q} \underset{\substack{d_1, d_1' \ldots, d_k, d_k' \geq 1\\q_w, W, [d_i,d_i'] \text{ pairwise coprime}}}{\sum \cdots \sum} \prod_{i=1}^k 
\frac{\mu(d_i) \mu(d_i')}{[d_i,d_i']} f_i \left( \frac{\log d_i}{\log q} \right)
g_i \left( \frac{\log d_i'}{\log q} \right),
\end{align*}
where $q_w:=\prod_{p|q, p > w} p$ and $W = \prod_{p \le w} p$. Then, applying \cite[Lemma 30]{MR3373710}, this becomes
\begin{align*}
T_1 &= (\mathfrak{I}(\boldsymbol{f},\boldsymbol{g})+o(1))\frac{\varphi(q)}{W_q} 
\left(  \frac{\varphi(W)}{W} \times   \frac{ \varphi(q_w)}{q_w}\times \log q \right)^{-k} \\
&= (\mathfrak{I}(\boldsymbol{f},\boldsymbol{g})+o(1))\frac{\varphi(q)}{W_q} 
\left(  \frac{\varphi(W_q)}{W_q} \times   \frac{ \varphi(q)}{q}\times \log q \right)^{-k},
\end{align*}
where $\boldsymbol{f}:=(f_1,\ldots,f_k), \boldsymbol{g}:=(g_1,\ldots,g_k)$ and
\begin{align*}
\mathfrak{I}(\boldsymbol{f},\boldsymbol{g})
:=\prod_{i=1}^k \int_{0}^{\infty}
f_i'(t_i) g_i'(t_i) dt_i.
\end{align*}

Recall that the supports of $(t_1, \ldots, t_k) \mapsto \prod_{i=1}^k f_{i}(t_i)$ and $(t_1, \ldots, t_k) \mapsto \prod_{i=1}^k g_{i}(t_i)$ are both contained in $\Delta_k(\theta, \delta)$. We bound $T_2$ using \eqref{eq:prop3.1} and \eqref{eq:estimate of S(d,d')}:
\begin{align*}
    T_2 &\ll_{f, g} \tau(q) \underset{d_1, d_1' \ldots, d_k, d_k' \geq 1}{\sum \cdots \sum}  1_{\sum_{i=1}^k\frac{\log d_i}{\log q} \le \frac{\theta}{2}(1-\frac{1}{2\log k})} \cdot  1_{\sum_{i=1}^k\frac{\log d'_i}{\log q} \le \frac{\theta}{2}(1-\frac{1}{2\log k})} \\
    &= \tau(q) \cdot  \Bigg(\sum_{a \le q^{\frac{\theta}{2}(1-\frac{1}{2\log k})}} \tau_k(a)\Bigg)^2 \\
    &\ll_{f, g, k} q^{\theta(1-\frac{1}{3\log k})} = o\left(\frac{\varphi(q)}{W_q} 
\left(  \frac{\varphi(W_q)}{W_q} \times   \frac{ \varphi(q)}{q}\times \log q \right)^{-k}\right).
\end{align*}
The proposition follows.
\end{proof}

\begin{proof}[Proof of Proposition \ref{prop:prime}]
Since $\rho$ is a minorant for the indicator of primes, we know that $\rho(a+qh_{i_0}, q) \le 0$ if $(a, q) > 1$. Therefore, 
\begin{equation}
\label{eq:remove coprimality}
    \sum_{\substack{q < a \le 2q\\(a,q)=1}} \rho(a+qh_{i_0}, q) w_a \ge \sum_{\substack{q < a \le 2q}} \rho(a+qh_{i_0}, q) w_a
\end{equation}

For our analysis, we suppose $i_0 = k$ (the cases $i_0 < k$ are treated in the same way). If $\rho(a+qh_k, q) \ne 0$, then any divisor $d | a+qh_k$ must either satisfy $d=1$ or $d\ge q^\xi$ by Definition \ref{def:set of minorants} \ref{item:supported on rough}. For all $j, j'$, since both $f_{k, j}$ and $f_{k, j'}$ are supported in $[0, \xi]$, we have $\lambda_{f_{k, j}}(a+qh_k)\lambda_{f_{k, j'}}(a+qh_k) = f_{k, j}(0)f_{k, j'}(0)$ whenever $\rho(a+qh_k, q) \ne 0$. Thus,
\begin{equation}
\label{eq:prop 3.2 red 1}
    \sum_{\substack{q < a \le 2q}} \rho(a+qh_{i_0}, q) w_a = \sum_{j=1}^J \sum_{j'=1}^J c_jc_{j'}f_{k, j}(0)f_{k, j'}(0) \widetilde T(j, j'),
\end{equation}
with 
\[
    \widetilde T(j, j') \coloneqq \sum_{\substack{q < a \le 2q \\ a \equiv b_0 \Mod{W_q}}} \rho(a+qh_k, q)\prod_{i=1}^{k-1}
\lambda_{f_{i,j}}(a+qh_i)\lambda_{f_{i,j'}}(a+qh_{i}).
\]

For general  $f_{i}, g_i \in \mathcal{C}^{\infty}_{c}([0,\infty))$ with $1 \leq i \leq k-1$, define 
\begin{align*} 
\widetilde{T} \coloneqq \sum_{\substack{q < a \le 2q \\ a \equiv b_0 \Mod{W_q}}} \rho(a+qh_k, q)\prod_{\substack{i=1}}^{k-1} 
\lambda_{f_i}(a+qh_i) \lambda_{g_i}(a+qh_i).
\end{align*}
If we prove the asymptotic formula
\begin{align} \label{eq:prop3.2}
    \widetilde T =  \frac{{\mathfrak{J}}_{k}(\boldsymbol{f}, \boldsymbol{g})+o(1)}{\varphi(W_q)} \left(\sum_{\substack{q < a \le 2q}} \rho(a+qh_k, q)\right)
\left(  \frac{\varphi(W_q)}{W_q} \times   \frac{\varphi(q)}{q} \times \log q \right)^{-(k-1)}
\end{align}
with 
\begin{align*}
{\mathfrak{J}}_{k}(\boldsymbol{f}, \boldsymbol{g}):=
\prod_{i=1}^{k-1} \int_{0}^{\infty}
f_i'(t_i) g_i'(t_i) \, dt_i,
\end{align*}
then, after recombining the finitely many tensor terms, the proposition follows from \eqref{eq:remove coprimality}, \eqref{eq:prop 3.2 red 1}, Property \ref{item:pnt for minorant} with $\ell = h_k + 1$ of Definition \ref{def:set of minorants}, and the identity
\[
    \sum_{j=1}^J \sum_{j'=1}^J c_jc_{j'} f_{k, j}(0) f_{k, j'}(0) \mathfrak{J}_{k}(\boldsymbol{f}_{j}, \boldsymbol{f}_{j'}) = \mathfrak{J}_{k}(\partial_{t_1 \cdots t_k} F)\ge 0.
\]

After interchanging the order of summation, we have
\begin{align*}
\widetilde T = \underset{d_1, d_1' \ldots, d_{k-1}, d'_{k-1} \geq 1 }{\sum \cdots \sum} \left( \prod_{i=1}^{k-1} 
\mu(d_i) \mu(d_i') f_i \left( \frac{\log d_i}{\log q} \right)
g_i \left( \frac{\log d_i'}{\log q} \right) \right) \widetilde{S}(\boldsymbol{d}, \boldsymbol{d'}),
\end{align*}
where 
\begin{align*}
\widetilde{S}(\boldsymbol{d}, \boldsymbol{d'}):=
\sum_{\substack{q < a \le 2q\\ [d_i,d_i'] \, | \, a+qh_i, 1 \leq i \leq k-1\\a \equiv b_0 \Mod{W_q}}}
\rho(a+qh_k, q).
\end{align*}
Arguing as in the proof of Proposition \ref{prop:nonprime}, there exists $c_{W_q, \boldsymbol{d}, \boldsymbol{d'}}$ such that 
\begin{equation}
    \label{eq:tilde S rewriting}
    \widetilde{S}(\boldsymbol{d}, \boldsymbol{d'})= \frac{1}{\varphi(\widetilde r_{W_q,\boldsymbol{d}, \boldsymbol{d'}})}\sum_{\substack{q < a \le 2q}} \rho(a+qh_k, q) + \Delta\Big(q, h_k+1, c_{W_q, \boldsymbol{d}, \boldsymbol{d'}} + qh_k, \widetilde r_{W_q, \boldsymbol{d}, \boldsymbol{d'}}\Big)
\end{equation}
by the Chinese Remainder Theorem, with $\widetilde r_{W_q, \boldsymbol{d}, \boldsymbol{d'}} \coloneqq W_q \prod_{i=1}^{k-1} [d_i, d_i']$ and 
\[
    \Delta(q, \ell, c, r) \coloneqq \sum_{\ell q < a \le (\ell+1)q} \rho(a, q)\left(1_{a \equiv c \Mod r} - \frac{1_{(a, r)=1}}{\varphi(r)}\right),
\]
if $q, W_q$ and all $[d_i,d_i']$ are pairwise coprime. Note that $a+qh_k$ and $\widetilde r_{W_q, \boldsymbol d, \boldsymbol{d'}}$ are coprime whenever $\rho(a+qh_k, q) \ne 0$ since all prime factors of $\widetilde r_{W_q, \boldsymbol{d}, \boldsymbol{d'}}$ are $\le q^{\xi}$, which is why we did not explicitly restrict the sum in \eqref{eq:tilde S rewriting} to $a$'s with $a+qh_k$ and $\widetilde r_{W_q, \boldsymbol{d}, \boldsymbol{d'}}$ being coprime. We define $\widetilde T_1$ as
\begin{align*}
    \widetilde T_1 &\coloneqq \frac{1}{\varphi(W_q)} \left(\sum_{q < a \le 2q} \rho(a+qh_k, q)\right) \\
    &\qquad \times\underset{\substack{d_1, d_1' \ldots, d_{k-1}, d_{k-1}' \geq 1 \\ q, W_q, [d_i,d_i'] \text{ pairwise coprime}}}{\sum \cdots \sum} \prod_{i=1}^{k-1} 
\frac{\mu(d_i) \mu(d_i')}{\varphi([d_i, d'_i])} f_i \left( \frac{\log d_i}{\log q} \right)
g_i \left( \frac{\log d_i'}{\log q} \right), 
\end{align*}
and $\widetilde T_2 \coloneqq \widetilde T - \widetilde T_1$. Note that 
\begin{align*}
\widetilde T_1 &= \frac{1}{\varphi(W_q)} \left(\sum_{q < a \le 2q} \rho(a+qh_k, q)\right) \\
    &\qquad \times\underset{\substack{d_1, d_1' \ldots, d_{k-1}, d_{k-1}' \geq 1 \\ q_w, W, [d_i,d_i'] \text{ pairwise coprime}}}{\sum \cdots \sum} \prod_{i=1}^{k-1} 
\frac{\mu(d_i) \mu(d_i')}{\varphi([d_i, d'_i])} f_i \left( \frac{\log d_i}{\log q} \right)
g_i \left( \frac{\log d_i'}{\log q} \right).
\end{align*}
where we recall from the proof of the previous proposition $q_w = \prod_{p|q, p > w} p$ and $W = \prod_{p \le w} p$. Put $B_q:=\frac{\varphi(W_q)}{W_q}\frac{\varphi(q)}{q}\log q$. We will prove that $\widetilde T_2=o\!\left(\frac{q}{\varphi(W_q)\log q}B_q^{-(k-1)}\right)$. Applying \cite[Lemma 30]{MR3373710}, we have
\begin{align*}
    \widetilde T_1 &=  \frac{{\mathfrak{J}}_{k}(\boldsymbol{f}, \boldsymbol{g})+o(1)}{\varphi(W_q)} \left(\sum_{\substack{q < a \le 2q}} \rho(a+qh_k, q)\right)
\left(  \frac{\varphi(W)}{W} \times   \frac{\varphi(q_w)}{q_w} \times \log q \right)^{-(k-1)} \\
&= \frac{{\mathfrak{J}}_{k}(\boldsymbol{f}, \boldsymbol{g})+o(1)}{\varphi(W_q)} \left(\sum_{\substack{q < a \le 2q}} \rho(a+qh_k, q)\right)
\left(  \frac{\varphi(W_q)}{W_q} \times   \frac{\varphi(q)}{q} \times \log q \right)^{-(k-1)}
\end{align*}

To bound $\widetilde T_2$, recall the product $P(z) \coloneqq \prod_{p \le z} p$. The product $\prod_{i=1}^{k-1}  f_i \left( \frac{\log d_i}{\log q} \right)
g_i \left( \frac{\log d_i'}{\log q} \right) $ vanishes unless $(\boldsymbol{d}, \boldsymbol{d'})$ belongs to the set 
\[
    \mathcal D = \Big\{(\boldsymbol{d}, \boldsymbol{d'}) \in \mathbb N^{k-1} \times \mathbb N^{k-1} : \widetilde r_{W_q,\boldsymbol{d}, \boldsymbol{d'}} \le q^{\theta-\epsilon_k} \text{ and } \widetilde r_{W_q,\boldsymbol{d}, \boldsymbol{d'}}|P(q^{\delta(\epsilon_k)})\Big\}
\]
for $\epsilon_k \coloneqq \frac{\theta}{4\log k}$. We have the upper bound
\[
    \widetilde T_2 \ll \underset{(\boldsymbol{d}, \boldsymbol{d'}) \in \mathcal D}{\sum} \ \bigg| \Delta\Big(q, h_k+1, c_{W_q, \boldsymbol{d}, \boldsymbol{d'}} + qh_k, \widetilde r_{W_q, \boldsymbol{d}, \boldsymbol{d'}}\Big) \bigg|.
\]
Let $\mathcal C$ be the set of positive integers $c \le P(q^{\delta(\epsilon_k)})$ satisfying $c\equiv b_0+qh_k \Mod{W_q}$ and for all $p \in (w, q^{\delta(\epsilon_k)}]$, we have $c \equiv q(h_k-h_i) \Mod p$ for some $i \in \{1, \ldots, k-1\}$. We have 
\[
    \widetilde T_2 \ll \sum_{(\boldsymbol{d}, \boldsymbol{d'}) \in \mathcal D} \max_{\substack{c \in \mathcal C \\ (c, \widetilde r_{W_q, \boldsymbol{d}, \boldsymbol{d'}}) = 1}} \bigg| \Delta\Big(q, h_k+1, c, \widetilde r_{W_q, \boldsymbol{d}, \boldsymbol{d'}}\Big) \bigg|
\]

Let $\omega(n)$ count the distinct prime factors of $n$. Note that for any squarefree $r\ge 1$ satisfying $W_q|r$ and $(r,q)=1$, we have
\begin{align*}
&\#\Big\{(\boldsymbol{d}, \boldsymbol{d'}) \in \mathbb N^{k-1} \times \mathbb N^{k-1} \,:\, r=\widetilde r_{W_q,\boldsymbol{d}, \boldsymbol{d'}}\Big\} \\
&\qquad = \sum_{a_1 \cdots a_{k-1} = r/W_q} \prod_{i=1}^{k-1}\#\Big\{(d_i, d'_i) \in \mathbb N \times \mathbb N \,:\, a_i=[d_i, d'_i]\Big\} = (3k-3)^{\omega(r/W_q)}.
\end{align*}

Therefore, 
\begin{align} \label{eq:prop3.2et}
\widetilde T_2 \ll \sum_{\substack{ r \leq q^{\theta-\epsilon_k} \\ r| P(q^{\delta(\epsilon_k)}) \\ (r,q)=1}} (3k-3)^{\omega(r)}
\max_{\substack{c \in \mathcal C \\ (c, r) = 1}} \bigg| \Delta\Big(q, h_k+1, c, r\Big) \bigg|.
\end{align}
For any given $r$ satisfying $r|P(q^{\delta(\epsilon_k)})$, $r \le q^{\theta - \epsilon_k}$, and $(r,q)=1$, and for any given $c_0 \in \mathcal C$ coprime to $r$, we have the lower bound
\[
    \#\{c \in \mathcal C : c \equiv c_0 \Mod{r}\} \ge \frac{|\mathcal C|}{k^{\omega(r)}}
\]
by the Chinese Remainder Theorem. Thus, 
\begin{align*}
\max_{\substack{c \in \mathcal C \\ (c, r) = 1}} \bigg| \Delta\Big(q, h_k+1, c, r\Big) \bigg| &\le \frac{k^{\omega(r)}}{|\mathcal C|} \sum_{\substack{c \in \mathcal C \\ (c, r) = 1}} \bigg|\Delta\Big(q, h_k+1, c, r\Big) \bigg| \\
&\ll \frac{q^{1/2}}{\varphi(r)^{1/2}}\cdot \frac{k^{\omega(r)}}{|\mathcal C|} \sum_{\substack{c \in \mathcal C \\ (c, r) = 1}} \bigg|\Delta\Big(q, h_k+1, c, r\Big) \bigg|^{1/2},
\end{align*}
where the last inequality follows from the trivial bound $|\Delta(q, h_k+1, c, r)| \ll q/\varphi(r)$ (we have this bound because of the property that $\rho$ is bounded above and below over its domain).

Combining this with \eqref{eq:prop3.2et} and using the Cauchy-Schwarz inequality, we have 
\begin{align*}
    \widetilde T_2 &\ll \frac{q^{1/2}}{|\mathcal C|} \sum_{c \in \mathcal C} \sum_{\substack{ r \leq q^{\theta-\epsilon_k} \\ r| P(q^{\delta(\epsilon_k)}) \\ (r,q)=1 \\ (c, r) = 1}} \frac{(3k(k-1))^{\omega(r)}}{\varphi(r)^{1/2}} \cdot \bigg|\Delta\Big(q, h_k+1, c, r\Big) \bigg|^{1/2} \\
    &\ll \frac{q^{1/2}}{|\mathcal C|} \left(\sum_{r \le q} \frac{(9k^4)^{\omega(r)}}{\varphi(r)}\right)^{1/2} \times \sum_{c \in \mathcal C} \Bigg(\sum_{\substack{ r \leq q^{\theta-\epsilon_k} \\ r| P(q^{\delta(\epsilon_k)}) \\ (r,q)=1 \\ (c, r) = 1}} \bigg|\Delta\Big(q, h_k+1, c, r\Big) \bigg|\Bigg)^{1/2} \ll \frac{q}{(\log q)^{k^4/2}}
\end{align*}
by applying Property \ref{item:level of distribution} of Definition \ref{def:set of minorants} with $A = 10k^4$. Since $\frac{q}{(\log q)^{k^4/2}}=o\!\left(\frac{q}{\varphi(W_q)\log q}B_q^{-(k-1)}\right)$, this proves the required error estimate and hence \eqref{eq:prop3.2}. The proposition follows.
\end{proof}

\begin{proof}[Proof of Proposition \ref{prop:small_prime}]
Similarly, it suffices to estimate
\begin{align} \label{eq:prop3.3}
  \sum_{\substack{q < a \le 2q \\ (a,q)=1 }} \Bigg(\sum_{\substack{p| a+qh_{k}  \\ p \leq q^{\gamma} }} 1 \Bigg) 1_{a \equiv b_0 \Mod{W_q}} 
   \prod_{i=1}^{k} 
\lambda_{f_i}(a+qh_i) \lambda_{g_i}(a+qh_i).
\end{align}
If $p\mid q$, then $(a,q)=1$ implies $p\nmid a+qh_k$, so these primes contribute nothing to \eqref{eq:prop3.3}. If $p\mid W_q$, the choice of $b_0$ likewise excludes divisibility; if $p\nmid qW_q$, the condition $p\mid a+qh_k$ selects one residue class modulo $p$. The unmarked local factors therefore match those in Proposition \ref{prop:nonprime}, including its $\varphi(q)/q$ normalization. Adapting the proof of \cite[Proposition 14 (Almost primality), pp.~23--27]{MR3373710}, one can show that
\begin{gather*}
\sum_{\substack{q < a \le 2q \\ (a,q)=1 }} 1_{p| a+qh_{k}} 1_{a \equiv b_0 \Mod{W_q}} 
\prod_{i=1}^{k} 
\lambda_{f_i}(a+qh_i) \lambda_{g_i}(a+qh_i) \\
\ll \frac{\log p}{p\log q} \times \frac{\varphi(q)}{W_q} 
\left(  \frac{\varphi(W_q)}{W_q} \times   \frac{ \varphi(q)}{q}\times \log q \right)^{-k}.
\end{gather*}
Then, by summing over primes $p,$ the expression (\ref{eq:prop3.3}) is
\begin{align*}
\ll \gamma \times \frac{\varphi(q)}{W_q} 
\left(  \frac{\varphi(W_q)}{W_q} \times   \frac{ \varphi(q)}{q}\times \log q \right)^{-k},
\end{align*}
and the proposition follows.
\end{proof}

\section{Proof of Theorem \ref{thm:dynamic}}
Let $h:=C\exp(3.8075m),$ where $C>0$ is the absolute constant in Corollary \ref{stadlman_corollary}. Then
using Lemma \ref{lem:poincare} followed by the triangle inequality, there exists the least integer $1 \leq q \leq \mu(B(x_0;\varepsilon/4h))^{-1}$ such that 
\begin{align*}
d(T^q x_0, x_0)<\frac{\varepsilon}{2h}.
\end{align*}
Since $x_0$ is nonperiodic, the least $q$ satisfying $d(T^q x_0,x_0)<\varepsilon/(2h)$ tends to infinity as $\varepsilon\to0$, so for $\varepsilon$ sufficiently small in terms of $m$ and $x_0$, Corollary \ref{stadlman_corollary} applies to this $q$.
Let 
$\mathcal{A}:=\{ 1 \leq a \leq q \,: \, p_m(q,a) \leq  q h, (a,q)=1 \}$ and $\mathcal{X}_{\mathcal{A}}:=\bigcup_{a \in \mathcal{A}} B(T^a x_0; \varepsilon/4h ).$ Note that $\mathcal{X}_{\mathcal{A}}$ is in fact a disjoint union by the minimality of $q.$ Then, for $a\in\mathcal{A}$, $x\in B(T^a x_0;\varepsilon/4h)$, and $1\leq i\leq m$, we have
\begin{align} \label{eq:ineq}
d(T^{p_i(q,a)} x_0, x) \leq & \ d(T^{p_i(q,a)} x_0, T^a x_0)
+ d(T^{a} x_0, x) \nonumber \\
< & \ d(T^{p_i(q,a)} x_0, T^a x_0)
+ \frac{\varepsilon}{4h}.
\end{align}
Since $T$ is an isometry, we have
\begin{align*}
d(T^{p_i(q,a)} x_0, T^a x_0) &=  d(T^{p_i(q,a)-a} x_0,  x_0) \\
& \leq \sum_{1 \leq j \leq \frac{p_i(q,a)-a}{q}} d(T^{jq}x_0,T^{(j-1)q}x_0) \\
& =   \frac{p_i(q,a)-a}{q} \times d(T^qx_0,x_0) \\
& < \frac{\varepsilon}{2}.
\end{align*}
Then, it follows from (\ref{eq:ineq}) that for all $x \in \mathcal X_\mathcal A$, there exists $a \in \mathcal A$ such that $d(T^{p_i(q,a)} x_0, x)< \varepsilon$ for all $i \in \{1, \ldots, m\}$, and thus 
$p_m^T(x_0,x;\varepsilon) \leq qh.$

On the other hand, since $T^{-a}(B(T^a x_0;\varepsilon/4h))=B(x_0;\varepsilon/4h)$ and $T$ is measure-preserving, we have 
$\mu(B(T^a x_0;\varepsilon/4h))=\mu(B(x_0;\varepsilon/4h)),$
which implies that
\begin{align*}
\mu(\mathcal{X}_{\mathcal{A}}) & = \# \mathcal{A} \times \mu(B(x_0;\varepsilon/4h)) \\
& \gg_{m} \varphi(q) (1+\log q)^{-C'\exp(3.8075m)} \times \mu(B(x_0;\varepsilon/4h))
\end{align*}
by Corollary \ref{stadlman_corollary}. Set $J=\lceil\log_2(4h)\rceil$; iterating pseudo-doubling gives $\mu(B(x_0;\varepsilon/4h))\ge\lambda^{-J}\mu(B(x_0;\varepsilon))$ and $q\le\lambda^J\mu(B(x_0;\varepsilon))^{-1}\ll_\lambda m^{\log_2\lambda}\exp(3.8075m\log_2\lambda)\mu(B(x_0;\varepsilon))^{-1}$. This proves the theorem.

To obtain Theorem \ref{thm:1.2} for any $x_0\in\operatorname{supp}\mu$, let $Q(m)$ be the modulus threshold in Corollary \ref{stadlman_corollary} and increase $h$ to
$h_m:=\max\{h,\max_{1\leq q<Q(m)}\min_{(a,q)=1}p_m(q,a)/q\}$, with the usual convention for $q=1$.
Choose the least return at radius $\varepsilon/(2h_m)$ as above. For $q\geq Q(m)$, Corollary \ref{stadlman_corollary} supplies a suitable residue class; for $q<Q(m)$, its existence follows from the definition of $h_m$. The same distance and pseudo-doubling estimates give Theorem \ref{thm:1.2} with an implied constant depending only on $m$ and $\lambda$.

\section{First return time} \label{sec:1streturn}

Given any $\varepsilon>0, x_0 \in X$ and $ x \in B(x_0 ; \varepsilon),$ let $\tau_{x_0;\varepsilon}(x) := \min\{n \geq 1: T^n x \in B(x_0 ; \varepsilon) \}$ denote the \textit{first return time}. For instance, the least integer $q \geq 1$ for which $d(T^q x_0, x_0)<\varepsilon/2h$ in the proof of Theorem \ref{thm:dynamic} is the first return time $\tau_{x_0;\varepsilon/2h}(x_0).$ Suppose the measure-preserving map $T$ is ergodic and $\mu(B(x_0;\varepsilon/2h))>0$. Then Kac's lemma \cite{MR0022323} states that
\begin{align*}
    \int_{B(x_0, \varepsilon/2h)} \tau_{x_0;\varepsilon/2h}(x) d\mu(x) = 1,
\end{align*}
i.e., under the conditional probability measure $\mu(\,\cdot\mid B(x_0;\varepsilon/2h))$, the expected return time is $\mu(B(x_0,\varepsilon/2h))^{-1}$. For the following conditional estimate, suppose instead that $x_0$ is nonperiodic and its first return time satisfies
\begin{align} \label{eq:1streturn}
\mu(B(x_0;\varepsilon/2h))^{-1}(\log \mu(B(x_0;\varepsilon/2h))^{-1} )^{-100}  \ll \tau_{x_0;\varepsilon/2h}(x_0) \ll \mu(B(x_0;\varepsilon/2h))^{-1}
\end{align}
This is a separate hypothesis on the centre $x_0$.

Suppose now (\ref{eq:1streturn}) holds. Let $q':=\tau_{x_0;\varepsilon/2h}(x_0)$. Then the proof of Theorem \ref{thm:dynamic} gives
\begin{gather*}
 \mu ( \{ x \in X \,: \, p_m^T(x_0,x;\varepsilon) \leq C {\mu(B(x_0;\varepsilon))}^{-1}  m^{\log_2 \lambda} \exp( 3.8075(1 + \log_2 \lambda)  m) \} )  \\
\gg_{m, \lambda} \varphi(q') (1+\log q')^{-C'\exp(3.8075m)} \mu(B(x_0;\varepsilon)) \\
\gg_{m,\lambda}  (\log \mu(B(x_0; \varepsilon))^{-1})^{-C''\exp(3.8075m)}
\end{gather*}
for some absolute constant $C''>0.$

Unfortunately, none of the circle rotations are even weakly mixing (see \cite[p. 51]{MR2723325} for instance). In fact, one can verify that the assumption (\ref{eq:1streturn}) fails for those circle rotations $T_{\alpha}:x \mapsto x+\alpha$ with $\alpha$ of \textit{type} $\eta_{\alpha}>1$,\footnote{By the Borel--Cantelli lemma, the set of irrationals $\alpha$ of type $\eta_{\alpha}>1$ has Lebesgue measure $0.$}
where
\begin{align*}
\eta_{\alpha}:=\sup \{ \theta > 0 \,: \, \liminf_{n \to \infty} n^{\theta} \|n\alpha\| =0 \}.
\end{align*}
More precisely, let $\tau_{\varepsilon}(\alpha):=\min\{n\geq 1 : \| n \alpha \| < \varepsilon\}.$ Then 
Choe--Seo \cite{MR1857291} proved that
$$ \liminf_{\varepsilon\to 0^{+}} \frac{\log \tau_{\varepsilon}(\alpha)}{\log \varepsilon^{-1}} = \eta_{\alpha}^{-1}. $$
Nevertheless, we are able to establish tight lower and upper bounds for the first return time for most $\alpha \in (0,1)$ of type $\eta_{\alpha} \leq 1,$ thereby strengthening \cite[Theorem 3]{MR1999584}.
\begin{proposition} \label{prop:1streturn}
Let $\alpha \in (0,1), \delta>0$ and $\varepsilon>0$ sufficiently small (in terms of $\alpha$ and $\delta$). Then
\begin{enumerate}
\item for almost all $\alpha \in (0,1)$ with respect to the Lebesgue measure, we have
\begin{align*}
\varepsilon^{-1}(\log \varepsilon^{-1})^{-(1+\delta)}  \ll_{\delta}
\tau_{\varepsilon}(\alpha)
\ll \varepsilon^{-1}; 
\end{align*}
\item the irrational $\alpha \in (0,1)$ has bounded partial quotients\footnote{An irrational $\alpha$ with continued fraction expansion $[a_0;a_1,a_2,\ldots]$ is said to have bounded partial quotients if $a_n$ is uniformly bounded for $n \geq 1.$ They are badly approximable by rationals, and vice versa, i.e., there exists a constant $c>0$ such that $|\alpha - \frac{p}{q}| > \frac{c}{q^2}$ for all $\frac{p}{q} \in \mathbb{Q}$. For instance, this includes quadratic irrationals but excludes the base of the natural logarithm, $e=[2;1,2,1,1,4,1,1,6,1,1,8,\dots].$ } if and only if
\begin{align*}
\tau_{\varepsilon}(\alpha) \asymp_{\alpha}  \varepsilon^{-1}.
\end{align*}
\end{enumerate}
\end{proposition}

In particular, the equation (\ref{eq:1streturn}) holds, allowing us to derive the following (metric) Diophantine approximation result on early visits at prime times.

\begin{corollary} \label{cor:almostall}
Let $m \geq 2$ and $\alpha \in (0,1).$ Then for any sufficiently small $\varepsilon>0$ (in terms of $m, \alpha$), there exist absolute constants $C, C'>0$ such that
\begin{enumerate}
\item for almost all $\alpha \in (0,1)$ with respect to the Lebesgue measure $\mu,$ we have
\begin{align*}
\mu ( \{ \beta \in [0,1] \,: \, p_m^{\alpha}(\beta;\varepsilon) \leq C \varepsilon^{-1} m \exp( 7.615 m) \} ) \gg_{m}  (\log \varepsilon^{-1})^{-C'\exp(3.8075m)};
\end{align*}
\item for $\alpha \in (0,1)$ with bounded partial quotients, we also have
\begin{align*}
\mu ( \{ \beta \in [0,1] \,: \, p_m^{\alpha}(\beta;\varepsilon) \leq C \varepsilon^{-1} m \exp(  7.615m) \} ) \gg_{m,\alpha}  (\log \varepsilon^{-1})^{-C'\exp(3.8075m)}.
\end{align*}
\end{enumerate}
\end{corollary}

\begin{proof}[Proof of Proposition \ref{prop:1streturn}]
We begin with the proof of the first part.
The upper bound follows from the Dirichlet approximation theorem. Since
$\sum_{n=2}^{\infty}1/(n(\log n)^{1+\delta})$ converges, the Borel--Cantelli lemma gives, for almost all irrational $\alpha \in [0,1]$ and all sufficiently large integers $n$,
\begin{align*}
\|n\alpha\|>\frac{1}{n(\log n)^{1+\delta}}.
\end{align*}
For such an $\alpha$, put $n=\tau_\varepsilon(\alpha)$. As $\varepsilon\to0$, we have $n\to\infty$, while Dirichlet's theorem gives $n\leq 2\varepsilon^{-1}$ for sufficiently small $\varepsilon$. Thus
\begin{align*}
\varepsilon>\|n\alpha\|>\frac{1}{n(\log n)^{1+\delta}}
\geq \frac{1}{n(\log(2\varepsilon^{-1}))^{1+\delta}},
\end{align*}
which gives the claimed lower bound.

If $\alpha$ has bounded partial quotients, then $\|n\alpha\|\gg_\alpha n^{-1}$ for every $n\geq1$. Taking $n=\tau_\varepsilon(\alpha)$ and using $\|n\alpha\|<\varepsilon$ gives $\tau_\varepsilon(\alpha)\gg_\alpha\varepsilon^{-1}$. The upper bound follows from Dirichlet's theorem.

For the converse, let $\alpha \in (0,1)$ be an irrational number with unbounded partial quotients and denote $q_{n_k}:=\min\{ q_n \,: \, q_{n+1} > kq_n \}$ for $k\geq 1.$ Note that such $q_{n_k}$ always exists. Take $\varepsilon_k:=q_{n_{k+1}}^{-1}$ for $k\geq 1.$ Then by definition
\begin{align*}
\tau_{\varepsilon_k}(\alpha) \leq q_{n_k} < \frac{q_{n_{k+1}}}{k} = \frac{\varepsilon_{k}^{-1}}{k}.
\end{align*}
Therefore, it is impossible that $\tau_{\varepsilon}(\alpha) \gg_{\alpha} \varepsilon^{-1}$ as $\varepsilon \to 0^{+},$ and the proof is completed.
\end{proof}

\section*{Acknowledgements}
The authors are grateful to Andrew Granville and Dimitris Koukoulopoulos for their advice and encouragement. They would also like to thank Nikos Frantzikinakis and Bryna Kra for their helpful comments. T. H. is supported by the
Research Council of Finland grants nos. 364214 and 370133. S.-K. L. is supported by the Croucher Fellowship for Postdoctoral Research.

\printbibliography

\end{document}